\documentclass{proc-l}

\makeatletter
\def\@printyear{TT}
\makeatother

\issueinfo{138}{08}{August}{2010}
\pagespan{2863}{2872}
\dateposted{March 23, 2010}
\PII{S 0002-9939(10)10284-6}
\copyrightinfo{2010}{American Mathematical Society}

\revertcopyright

\usepackage{url}

\newcommand{\R}{\mathbb{R}}

\newcommand{\co}{\colon}
\newcommand{\vol}{\mathop{\textup{vol}}}
\newcommand{\Aut}{\mathop{\textup{Aut}}}
\newcommand{\supp}{\mathop{\textup{supp}}}

\theoremstyle{plain}
\newtheorem{theorem}{Theorem}

\newtheorem{lemma}[theorem]{Lemma}

\newtheorem{conjecture}[theorem]{Conjecture}

\theoremstyle{definition}

\newtheorem*{ex}{Examples}

\numberwithin{theorem}{section}
\numberwithin{equation}{section}

\begin{document}
\title[Asymmetric yet balanced point configurations]
{Point configurations\linebreak[1]  that are
asymmetric yet balanced}

\author[H. Cohn]{Henry Cohn}
\address{Microsoft Research New England,
One Memorial Drive,
Cambridge, Massachusetts 02142}
\email{cohn@microsoft.com}

\author[N. D. Elkies]{Noam D. Elkies}
\address{Department of Mathematics,
Harvard University,
Cambridge, Massachusetts 02138}
\email{elkies@math.harvard.edu}

\author[A. Kumar]{Abhinav Kumar}
\address{Department of Mathematics,
Massachusetts Institute of Technology,
Cambridge, Massachusetts 02139}
\email{abhinav@math.mit.edu}

\author[A. Sch\"urmann]{Achill Sch\"urmann}
\address{Institute of Applied Mathematics,
TU Delft,
Mekelweg 4,
2628 CD Delft, The Netherlands}
\email{a.schurmann@tudelft.nl}
\thanks{The first, second, and third authors
 thank the Aspen Center for Physics for its
hospitality and support.  The first, third, and fourth authors
 thank the Hausdorff Research
Institute for Mathematics. The third and fourth authors
 thank Microsoft Research.
The third author
 was supported in part by a Clay Liftoff Fellowship and by National
Science Foundation grant DMS-0757765.  The second author  was supported
in part by NSF grant DMS-0501029, and the fourth author  was
supported in part by Deutsche Forschungsgemeinschaft grant SCHU
1503/4-2.}

\subjclass[2010]{Primary 52B15; Secondary 05B40, 52C17, 82B05}
\date{March 30, 2009}
\commby{Jim Haglund}
\begin{abstract}
A configuration of particles confined to a sphere is balanced if it is
in equilibrium under all force laws (that act between pairs of points
with strength given by a fixed function of distance). It is
straightforward to show that every sufficiently symmetrical
configuration is balanced, but the converse is far from obvious. In
1957 Leech completely classified the balanced configurations in $\mathbb{R}^3$,
and his classification is equivalent to the converse for~$\mathbb{R}^3$. In
this paper we disprove the converse in high dimensions. We construct
several counterexamples, including one with trivial symmetry group.
\end{abstract}
\maketitle

\section{Introduction}

A finite configuration of points on the unit sphere $S^{n-1}$ in $\R^n$
is \emph{balanced} if it is in equilibrium (possibly unstable) under
all pairwise forces depending only on distance, assuming the points are
confined to the surface of the sphere. In other words, the net forces
acting on the points are all orthogonal to the sphere.  As is usual in
physics, any two distinct particles exert forces on each other,
directed oppositely and with magnitude equal to some fixed function of
the Euclidean distance between them. The net force on each point is the
sum of the contributions from the other points.

For example, the vertices of any regular polyhedron are balanced.  On
the other hand, most configurations are not balanced.  Even if some
points are in equilibrium under one force law, there is no reason to
expect that they will be in equilibrium under every force law, and
usually they will not be.  The balanced configurations are quite
remarkable.

The condition of being balanced was defined by Leech in \cite{Leech}.
It arises in the search for energy-minimizing point configurations on
spheres.  Given a potential function, typically an inverse-power law,
how should we arrange some particles to minimize the total potential
energy? This problem originated in Thomson's model of the atom in
\cite[p.~255]{T}. Of course, that model was superseded by quantum
mechanics, but it remains of considerable mathematical interest. It
provides a natural measure of how well distributed points are on the
surface of the sphere, and it also offers the possibility of
characterizing important or beautiful configurations via extremal
properties.

In most cases the optimal configuration depends on the potential
function, but occasionally it does not. In \cite{CK}, Cohn and Kumar
introduced the concept of \emph{universally optimal} configurations,
which minimize energy not only under all inverse-power laws but also
under the broader class of completely monotonic potential functions (as
functions of squared Euclidean distance). In $\R^2$ the vertices of any
regular polygon form a universally optimal configuration. The vertex
sets of the regular tetrahedron, octahedron, or icosahedron are
universally optimal, but there are no larger examples in $\R^3$.
Higher-dimensional examples include the vertices of the regular simplex
and cross polytope (or hyperoctahedron), and also various exceptional
examples, notably the vertices of the regular $600$-cell, the $E_8$
root system, the Schl\"afli configuration of $27$ points in $\R^6$
corresponding to the $27$ lines on a cubic surface, and the minimal
vectors of the Leech lattice. A number of the sporadic finite simple
groups act on universal optima. See Tables~1 and~2 in \cite{BBCGKS} for
a list of the known and conjectured universal optima, as well as for a
discussion of how many more there might be. (They appear to be quite
rare.)

Every universal optimum is balanced (as we will explain below), but
balanced configurations do not necessarily minimize energy even
locally.  In the space of configurations, balanced configurations are
universal critical points for energy, but they are frequently saddle
points.  For example, the vertices of a cube are balanced, but one can
lower the energy by rotating the vertices of a facet. Nevertheless,
being balanced is an important necessary condition for universal
optimality.

The simplest reason a configuration would be balanced is due to its
symmetry: the net forces inherit this symmetry, which can constrain
them to point orthogonally to the sphere. More precisely, call a finite
subset $\mathcal{C} \subset S^{n-1}$ \emph{group-balanced} if for every
$x \in \mathcal{C}$, the stabilizer of $x$ in the isometry group of
$\mathcal{C}$ fixes no vectors in $\R^n$ other than the multiples of
$x$. A group-balanced configuration must be balanced, because the net
force on $x$ is invariant under the stabilizer of $x$ and is thus
orthogonal to the sphere.

In his 1957 paper \cite{Leech}, Leech completely classified the
balanced configurations in $S^2$.  His classification shows that they
are all group-balanced, and in fact the complete list can be derived
easily from this assertion using the classification of finite subgroups
of $O(3)$. However, Leech's proof is based on extensive case analysis,
and it does not separate cleanly in this way. Furthermore, the
techniques do not seem to apply to higher dimensions.

It is natural to wonder whether all balanced configurations are
group-balanced in higher dimensions.  If true, that could help explain
the symmetry of the known universal optima.  However, in this paper we
show that balanced configurations need not be group-balanced. Among
several counterexamples, we construct a configuration of $25$ points in
$\R^{12}$ that is balanced yet has no nontrivial symmetries.

This result is compatible with the general philosophy that it is
difficult to find conditions that imply symmetry in high dimensions,
short of simply imposing the symmetry by fiat.  We prove that if a
configuration is a sufficiently strong spherical design relative to
the number of distances between points in it, then it is automatically
balanced (see Theorem~\ref{theorem:main}).  Every spectral embedding of
a strongly regular graph satisfies this condition (see
Section~\ref{section:srg}). There exist strongly regular graphs with no
nontrivial symmetries, and their spectral embeddings are balanced but
not group-balanced.

Before we proceed to the proofs, it is useful to rephrase the condition
of being balanced as follows: a configuration $\mathcal{C}$ is balanced
if and only if for every $x \in \mathcal{C}$ and every real number $u$,
the sum $S_u(x)$ of all $y \in \mathcal{C}$ whose inner product with
$x$ is $u$ is a multiple of~$x$. The reason is that the contribution to
the net force on~$x$ from the particles at a fixed distance is in the
span of $x$ and $S_u(x)$. Since we are using arbitrary force laws, each
contribution from a fixed distance must independently be orthogonal to
the sphere (since we can weight them however we desire). Note that a
group-balanced configuration $\mathcal{C}$ clearly satisfies this
criterion: for every $x \in \mathcal{C}$ and every real number $u$, the
sum $S_u(x)$ is itself fixed by the stabilizer of $x$ and hence must be
a multiple of $x$.

An immediate consequence of this characterization of balanced
configurations is that it is easy to check whether a given
configuration is balanced. By contrast, it seems difficult to check
whether a configuration is universally optimal.  For example, the paper
\cite{BBCGKS} describes a $40$-point configuration in $\R^{10}$ that
appears to be universally optimal, but so far no proof is known.

So far we have not explained why universal optima must be balanced. Any
optimal configuration must be in equilibrium under the force laws
corresponding to the potential functions it minimizes, but no
configuration could possibly minimize all potential functions
simultaneously (universal optima minimize a large but still restricted
class of potentials). The explanation is that a configuration is
balanced if and only if it is balanced for merely the class of
inverse-power force laws. In the latter case, we cannot weight the
force contributions from different distances independently. However, as
the exponent of the force law tends to infinity, the force contribution
from the shortest distance will dominate unless it acts orthogonally to
the sphere. This observation can be used to isolate each force
contribution in order by distance.  Alternatively, we can argue that
the configuration is balanced under any linear combination of
inverse-power laws and hence any polynomial in the reciprocal of
distance.  We can then isolate any single distance by choosing that
polynomial to vanish at all the other distances.

\section{Spherical designs}

Recall that a \emph{spherical $t$-design} in $S^{n-1}$ is a (nonempty)
finite subset $\mathcal{C}$ of $S^{n-1}$ such that for every polynomial
$p \co \R^n \to \R$ of total degree at most~$t$, the average of $p$
over $\mathcal{C}$ equals its average over all of $S^{n-1}$. In other
words,
$$
\frac{1}{|\mathcal{C}|}\sum_{x \in \mathcal{C}} p(x)
= \frac{1}{\vol(S^{n-1})} \int_{S^{n-1}} p(x) \, d\mu(x),
$$
where $\mu$ denotes the surface measure on $S^{n-1}$ and
$\vol(S^{n-1})$ is of course not the volume of the enclosed ball but
rather $\int_{S^{n-1}} d\mu(x)$.

\begin{theorem} \label{theorem:main}
Let $\mathcal{C} \subset S^{n-1}$ be a spherical $t$-design.  If for
each $x \in \mathcal{C}$, $$|\{\langle x,y \rangle : y \in \mathcal{C},
y \ne \pm x\}| \le t,$$ then $\mathcal{C}$ is balanced.
\end{theorem}

Here, $\langle\cdot,\cdot\rangle$ denotes the usual Euclidean inner
product.

\begin{proof}
Let $x$ be any element of $\mathcal{C}$, and let $u_1,\dots,u_k$ be the
inner products between $x$ and the elements of $\mathcal{C}$ other than
$\pm x$.  By assumption, $k \le t$.  We wish to show that for each $i$,
the sum $S_{u_i}(x)$ of all $z \in \mathcal{C}$ such that $\langle z,x
\rangle = u_i$ is a multiple of $x$.

Given any vector $y \in \R^n$ and integer $i$ satisfying $1 \le i \le
k$, define the degree $k$ polynomial $p \co \R^n \to \R$ by
$$
p(z) = \langle y,z \rangle \prod_{j \,:\, 1 \le j \le k, \, j \ne i}
 \big(\langle x,z \rangle - u_j\big).
$$

Suppose now that $y$ is orthogonal to $x$. Then the average of $p$ over
$S^{n-1}$ vanishes, because on the cross sections of the sphere on
which $\langle x,z \rangle$ is constant, each factor $\langle x,z
\rangle - u_j$ is constant, while $\langle y,z \rangle$ is an odd
function on such a cross section. More precisely, under the map $z
\mapsto 2 \langle x,z \rangle x - z$ (which preserves the component
of~$z$ in the direction of~$x$ and multiplies everything orthogonal to
$x$ by $-1$), the inner product with~$x$ is preserved while the inner
product with~$y$ is multiplied by~$-1$. Since $\mathcal{C}$ is a
$t$-design, it follows that the sum of $p(z)$ over $z \in \mathcal{C}$
also vanishes.

Most of the terms in this sum vanish: when $z = \pm x$, we have
$\langle y,z \rangle = 0$, and when $\langle x,z \rangle = u_j$, the
product vanishes unless $j=i$. It follows that the sum of $p(z)$ over
$z \in \mathcal{C}$ equals
$$
\sum_{z \in \mathcal{C} \,:\, \langle z,x \rangle = u_i}
\langle y,z \rangle
\prod_{j \,:\, 1 \le j \le k, \, j \ne i} (u_i - u_j)
=
\left(\prod_{j \,:\, 1 \le j \le k, \, j \ne i} (u_i - u_j)\right)
\big\langle y, S_{u_i}(x) \big\rangle.
$$

Because the first factor is nonzero, we conclude that $S_{u_i}(x)$ must
be orthogonal to $y$.  Because this holds for all $y$ orthogonal
to~$x$, it follows that $S_{u_i}(x)$ is a multiple of $x$, as desired.
\end{proof}

\begin{ex}
The vertices of a cube form a spherical $3$-design,
and only two inner products other than $\pm1$ occur between them, so
Theorem~\ref{theorem:main} implies that the cube is balanced. On the
other hand, not every group-balanced configuration satisfies the
hypotheses of the theorem.  For example, the configuration in $S^2$
formed by the north and south poles and a ring of $k$ equally spaced
points around the equator is group-balanced, but it is not even a
$2$-design if $k \ne 4$.  In Section~\ref{section:srg} we will show
that Theorem~\ref{theorem:main} applies to some configurations that are
not group-balanced, so the two sufficient conditions for being balanced
are incomparable.
\end{ex}

\section{Counterexamples from strongly regular graphs}
\label{section:srg}

Every spectral embedding of a strongly regular graph is both a
spherical $2$-design and a $2$-distance set, so by
Theorem~\ref{theorem:main} they are all balanced.  Recall that to form
a spectral embedding of a strongly regular graph with $N$ vertices, one
orthogonally projects the standard orthonormal basis of $\R^N$ to a
nontrivial eigenspace of the adjacency matrix of the graph.  See
Sections~2 and~3 of \cite{CGS} for a brief review of the theory of
spectral embeddings. Theorem~4.2 in \cite{CGS} gives the details of the
result that every spectral embedding is a $2$-design, a fact previously
noted as part of Example~9.1 in \cite{DGS}.

\begin{figure}
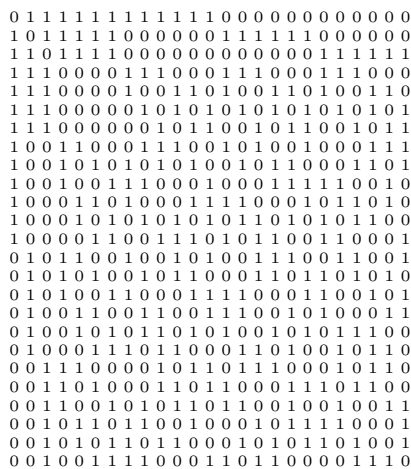

{\tiny \begin{center}
\begin{tabular}{ccccccccccccccccccccccccc}
0 1 1 1 1 1 1 1 1 1 1 1 1 0 0 0 0 0 0 0 0 0 0 0 0\\
1 0 1 1 1 1 1 0 0 0 0 0 0 1 1 1 1 1 1 0 0 0 0 0 0\\
1 1 0 1 1 1 1 0 0 0 0 0 0 0 0 0 0 0 0 1 1 1 1 1 1\\
1 1 1 0 0 0 0 1 1 1 0 0 0 1 1 1 0 0 0 1 1 1 0 0 0\\
1 1 1 0 0 0 0 1 0 0 1 1 0 1 0 0 1 1 0 1 0 0 1 1 0\\
1 1 1 0 0 0 0 0 1 0 1 0 1 0 1 0 1 0 1 0 1 0 1 0 1\\
1 1 1 0 0 0 0 0 0 1 0 1 1 0 0 1 0 1 1 0 0 1 0 1 1\\
1 0 0 1 1 0 0 0 1 1 1 0 0 1 0 1 0 0 1 0 0 0 1 1 1\\
1 0 0 1 0 1 0 1 0 1 0 1 0 0 1 0 1 1 0 0 0 1 1 0 1\\
1 0 0 1 0 0 1 1 1 0 0 0 1 0 0 0 1 1 1 1 1 0 0 1 0\\
1 0 0 0 1 1 0 1 0 0 0 1 1 1 1 0 0 0 1 0 1 1 0 1 0\\
1 0 0 0 1 0 1 0 1 0 1 0 1 0 1 1 0 1 0 1 0 1 1 0 0\\
1 0 0 0 0 1 1 0 0 1 1 1 0 1 0 1 1 0 0 1 1 0 0 0 1\\
0 1 0 1 1 0 0 1 0 0 1 0 1 0 0 1 1 1 0 0 1 1 0 0 1\\
0 1 0 1 0 1 0 0 1 0 1 1 0 0 0 1 1 0 1 1 0 1 0 1 0\\
0 1 0 1 0 0 1 1 0 0 0 1 1 1 1 0 0 0 1 1 0 0 1 0 1\\
0 1 0 0 1 1 0 0 1 1 0 0 1 1 1 0 0 1 0 1 0 0 0 1 1\\
0 1 0 0 1 0 1 0 1 1 0 1 0 1 0 0 1 0 1 0 1 1 1 0 0\\
0 1 0 0 0 1 1 1 0 1 1 0 0 0 1 1 0 1 0 0 1 0 1 1 0\\
0 0 1 1 1 0 0 0 0 1 0 1 1 0 1 1 1 0 0 0 1 0 1 1 0\\
0 0 1 1 0 1 0 0 0 1 1 0 1 1 0 0 0 1 1 1 0 1 1 0 0\\
0 0 1 1 0 0 1 0 1 0 1 1 0 1 1 0 0 1 0 0 1 0 0 1 1\\
0 0 1 0 1 1 0 1 1 0 0 1 0 0 0 1 0 1 1 1 1 0 0 0 1\\
0 0 1 0 1 0 1 1 0 1 1 0 0 0 1 0 1 0 1 1 0 1 0 0 1\\
0 0 1 0 0 1 1 1 1 0 0 0 1 1 0 1 1 0 0 0 0 1 1 1 0
\end{tabular}
\end{center}}
\caption{An adjacency matrix of a (25,12,5,6) strongly regular graph
with no nontrivial automorphisms.}\label{figure:srg25}
\end{figure}

The symmetry group of such a configuration is simply the combinatorial
automorphism group of the graph, so it suffices to find a strongly
regular graph with no nontrivial automorphisms. According to Brouwer's
tables \cite{B1}, the smallest such graph is a $25$-vertex graph with
parameters $(25,12,5,6)$ (the same as those of the Paley graph for the
\hbox{$25$-element} field), which has a spectral embedding in
$\R^{12}$. See Figure~\ref{figure:srg25} for an adjacency matrix.
Verifying that this graph has no automorphisms takes a moderate amount
of calculation, which is  best done by computer.

In fact, there are two such graphs with no nontrivial automorphisms
(the other is the complement of the graph in
Figure~\ref{figure:srg25}). Paulus classified the $(25,12,5,6)$
strongly regular graphs in \cite{P}; unfortunately, his paper was never
published. There are fifteen such graphs, whose automorphism groups
have a variety of sizes: two have order $1$, four have order $2$, two
have order $3$, four have order $6$, two have order $72$, and one has
order $600$ (the Paley graph). See \cite{B2} for more information.

The Paulus graphs give the lowest-dimensional balanced configurations
we have found that have trivial symmetry groups.  However, there are
lower-dimensional counterexamples (with some symmetry but not enough to
be group-balanced).  The lowest-dimensional one we have constructed is
in $\R^7$, and it can be built as follows; fortunately, no computer
calculations are needed.

Let $\mathcal{C}_n$ consist of the $n(n+1)/2$ midpoints of the edges of
a regular simplex in $\R^n$ (scaled so that $\mathcal{C}_n \subset
S^{n-1}$).  This configuration is a $2$-distance set, with the
distances corresponding to whether the associated edges of the simplex
intersect or not.  To compute the inner products, note that if
$x_1,\dots,x_{n+1}$ are the vertices of a regular simplex with
$|x_i|^2=1$ for all $i$ (and hence $\langle x_i,x_j\rangle=-1/n$ for $i
\ne j$), then for $i \ne j$ and $k \ne \ell$,
$$
\big\langle x_i+x_j,x_k+x_\ell \big\rangle = \begin{cases}
2-2/n & \textup{if $\{i,j\} = \{k,\ell\}$,}\\
1-3/n & \textup{if $|\{i,j\} \cap \{k,\ell\}|=1$, and}\\\
-4/n & \textup{if $\{i,j\} \cap \{k,\ell\} = \emptyset$.}
\end{cases}
$$
Thus, when we renormalize the vectors $x_i+x_j$ to lie on the unit
sphere, we find that the inner products between them are
$(1-3/n)/(2-2/n) = (n-3)/(2n-2)$ and $-(4/n)/(2-2/n) = -2/(n-1)$.

For $n>3$, the symmetry group of $\mathcal{C}_n$ is the same as that of
the original simplex (namely the symmetric group on the vertices of the
simplex).  Clearly, that group acts on $\mathcal{C}_n$. To see that
$\mathcal{C}_n$ has no other symmetries, we will show that the original
simplex can be constructed from it in such a way as to be preserved by
all symmetries of $\mathcal{C}_n$. Specifically, consider the subsets
of $\mathcal{C}_n$ of size $n$ in which all pairs of points are at the
minimal distance; the sums of these subsets are proportional to the
vertices of the original simplex. To see why, note that such a subset
corresponds to a collection of $n$ pairwise-intersecting edges of the
original simplex.  They must be exactly the edges containing one of the
vertices of the simplex: once two intersecting edges are specified,
only one other edge can intersect both without containing their common
vertex, so at most three edges can intersect pairwise without
containing a common vertex. (Note that this conclusion genuinely
requires that $n>3$, because $\mathcal{C}_3$ is an octahedron, which
has more symmetry than the tetrahedron from which it was derived.)

When $n=7$ the inner products in $\mathcal{C}_7$ are simply $\pm 1/3$.
The coincidence that these inner products are negatives of each other
is deeper than it appears, and it plays a role in several useful
constructions. For example, the union of $\mathcal{C}_7$ and its
antipode $-\mathcal{C}_7$ is a $3$-distance set, while in other
dimensions it would be a $5$-distance set.  In fact, $\mathcal{C}_7
\cup (-\mathcal{C}_7)$ is the unique $56$-point universal optimum in
$\R^7$, and it is invariant under the Weyl group of~$E_7$.  We will
make use of the unusual inner products in $\mathcal{C}_7$ to construct
a modification of it that is balanced but not group-balanced.

Within $\mathcal{C}_7$, there are regular tetrahedra (i.e., quadruples
of points with all inner products $-1/3$). Geometrically, such a
tetrahedron corresponds to a set of four disjoint edges in the original
simplex, and there is a unique such set up to symmetry, since the
simplex in $\R^7$ has eight vertices and all permutations of these
vertices are symmetries. Choose a set of four disjoint edges and call
them the distinguished edges.

We now define a modified configuration $\mathcal{C}_7'$ by replacing
each point in this tetrahedron by its antipode. Replacing the regular
tetrahedron preserves the $2$-design property, because the tetrahedron
is itself a $2$-design within the $2$-sphere it spans. In particular,
for every polynomial of total degree at most $2$, its sum over the
original tetrahedron is the same as its sum over the antipodal
tetrahedron.  Furthermore, when we replace the tetrahedron, all inner
products remain $\pm 1/3$ (some are simply multiplied by $-1$). Thus,
the resulting configuration $\mathcal{C}'_7$ remains both a
$2$-distance set and a $2$-design, so it is balanced by
Theorem~\ref{theorem:main}.

However, the process of inverting a tetrahedron reduces the symmetry
group.

\begin{lemma}
\label{lemma:group} The configuration $\mathcal{C}'_7$ has only
$4!\cdot 2^4 = 384$ symmetries, namely the permutations of the vertices
of that original simplex that preserve the set of four distinguished
edges.
\end{lemma}

\begin{proof}
There are clearly $4! \cdot 2^4$ symmetries of $\mathcal{C}_7$ that
preserve the set of distinguished edges of the simplex (they can be
permuted arbitrarily and their endpoints can be swapped).  All of these
symmetries preserve $\mathcal{C}'_7$.

To show that there are no further symmetries, it suffices to show that
the distinguished tetrahedron in $\mathcal{C}'_7$ is preserved under
all symmetries.  (For then the antipodal tetrahedron is also preserved,
and hence $\mathcal{C}_7$ is preserved as well.) Label the vertices of
the original simplex $1,2,\dots,8$, and suppose that the distinguished
edges correspond to the pairs $12$, $34$, $56$, and $78$. Label the
points of $\mathcal{C}'_7$ by the pairs for the corresponding edges.

There are at most two orbits under the symmetry group of
$\mathcal{C}'_7$, one containing $12$, $34$, $56$, and $78$ and the
other containing the remaining points.  We wish to show that these sets
do not in fact form a single orbit.  To separate the two orbits, we
will count the number of regular tetrahedra each point is contained in.
(We drop the word ``regular'' below.)  The answer will be seven for the
four distinguished points and eleven for the other points, so they
cannot lie in the same orbit.

Before beginning, we need a criterion for when the inner product
between two points in $\mathcal{C}'_7$ is $-1/3$.  If both points are
distinguished or both are nondistinguished, then that occurs exactly
when their label pairs are disjoint.  If one point is distinguished and
the other is not, then it occurs exactly when their label pairs
intersect.

Now it is straightforward to count the tetrahedra containing a
distinguished point; without loss of generality, call the point
 $12$.  There is one
tetrahedron of distinguished points, namely $\{12,34,56,78\}$.  If we
include a second distinguished point, say $34$, then there are two ways
to complete the tetrahedron using two nondistinguished points, namely
$\{12,34,13,24\}$ and $\{12,34,14,23\}$ (the two additional pairs must
be disjoint and each intersect both $12$ and $34$).  Because there are
three choices for the second distinguished point, this yields six
tetrahedra.  Finally, it is impossible to form a tetrahedron using
point $12$
and three nondistinguished points (one cannot choose three disjoint
pairs that each intersect $12$).  Thus, $12$ is contained in seven
tetrahedra.

To complete the proof, we need only show that a nondistinguished
point, without loss of generality $13$, is contained in more than seven
tetrahedra. There is a unique tetrahedron containing point $13$ and two
distinguished points, namely $\{13, 12, 34, 24\}$.  (There are only two
distinguished points that overlap with $13$, namely $12$ and $34$; then
the fourth point $24$ is determined.)  No tetrahedron can contain a
single distinguished point, as we saw in the previous paragraph, and if
a tetrahedron contains three distinguished points, then it must contain
the fourth. Thus, the only remaining possibility is that all the points
are nondistinguished.  The three other points in the tetrahedron must
be labeled with disjoint pairs from $\{2,4,5,6,7,8\}$, and the labels
$56$ and $78$ are not allowed (because those points are distinguished).
There are $6!/(2!^3\cdot3!)=15$ ways to split $\{2,4,5,6,7,8\}$ into
three disjoint pairs.  Among them, three contain the pair $56$, three
contain the pair $78$, and one contains both pairs.  Thus, there are
$15-3-3+1=10$ possibilities containing neither $56$ nor~$78$.  In
total, the point $13$ is contained in eleven tetrahedra.  Since it is
contained in more than seven tetrahedra, we see that points $12$ and $13$ are
in different orbits, as desired.
\end{proof}

By Lemma~\ref{lemma:group}, there are two orbits of points in
$\mathcal{C}'_7$, namely the four points in the tetrahedron and the
remaining $24$ points. The stabilizer of any point in the large orbit
actually fixes two such points. Specifically, consider the edge in the
original simplex that corresponds to the point. It shares its vertices
with two of the four distinguished edges (each vertex is in a unique
distinguished edge), and there is another edge that connects the other
two vertices of those distinguished edges.  For example, in the
notation of the proof of Lemma~\ref{lemma:group}, the edge $13$ has the
companion~$24$. This second edge has the same stabilizer as the first.
It follows that $\mathcal{C}'_7$ is not group-balanced.

If we interpret the configuration $\mathcal{C}'_7$ as a graph by using
its shorter distance to define edges, then we get a strongly regular
graph, with parameters $(28,12,6,4)$, the same as those of
$\mathcal{C}_7$. In fact, every $2$-design $2$-distance set yields a
strongly regular graph by Theorem~7.4 of \cite{DGS}. We have checked
using Brouwer's list \cite{B1} that spectral embeddings of strongly
regular graphs do not yield counterexamples in lower dimensions.  It
suffices to consider graphs with at most $27$ vertices, since by
Theorem~4.8 in \cite{DGS} no two-distance set in $S^5$ contains more
than $27$ points. Aside from the degenerate case of complete
multipartite graphs and their complements, the full list of strongly
regular graphs with spectral embeddings in six or fewer dimensions is
the pentagon, the Paley graph on $9$ vertices, the Petersen graph, the
Paley graph on $13$ vertices, the line graph of $K_6$, the Clebsch
graph, the Shrikhande graph, the $4 \times 4$ lattice graph, the line
graph of $K_7$, the Schl\"afli graph, and the complements of these
graphs.  It is straightforward to check that these graphs all have
group-balanced spectral embeddings. Of course there may be
low-dimensional counterexamples of other forms.

We suspect that there are no counterexamples in $\R^4$:

\begin{conjecture}
In $\R^4$, every balanced configuration is group-balanced.
\end{conjecture}

If true, this conjecture would lead to a complete classification of
balanced configurations in $\R^4$, because all the finite subgroups of
$O(4)$ are known (see for example \cite{CS}).  It is likely that, using
such a classification, one could prove completeness for the list of
known universal optima in $\R^4$, namely the regular simplices, cross
polytope, and $600$-cell, but we have not completed this calculation.

In $\R^5$ or $\R^6$, we are not willing to hazard a guess as to whether
all balanced configurations are group-balanced.  The construction of
$\mathcal{C}'_7$ uses such an ad hoc approach that it provides little
guidance about lower dimensions.

\section{Counterexamples from lattices}

In higher dimensions, we can use lattices to construct counterexamples
that do not arise from strongly regular graphs.  For example, consider
the lattice $\Lambda(G_2)$ in the Koch-\kern-.15exVenkov list of extremal even
unimodular lattices in $\R^{32}$ (see \cite[p.~212]{KV} or the
Nebe-Sloane catalogue~\cite{NS} of lattices).  This lattice has
$146880$ minimal vectors. When they are renormalized to be unit
vectors, only five inner products occur besides $\pm 1$ (namely, $\pm
1/2$, $\pm 1/4$, and $0$). By Corollary~3.1 of \cite{BV}, this
configuration is a spherical $7$-design. Hence, by
Theorem~\ref{theorem:main} it is balanced.  However,
$\Aut(\Lambda(G_2))$ is a relatively small group, of order $3 \cdot
2^{12}$, and one can check by computer calculations that some minimal
vectors have trivial stabilizers. (The lattice is generated by its
minimal vectors, and thus it and its kissing configuration have the
same symmetry group.) The kissing configuration of $\Lambda(G_2)$ is
therefore balanced but not group-balanced.

The case of $\Lambda(G_2)$ is particularly simple since some
stabilizers are trivial, but one can also construct lower-dimensional
counterexamples using lattices.  For example, let $L$ be the unique
$2$-modular lattice in dimension~$20$ with Gram matrix
determinant~$2^{10}$, minimal norm~$4$, and automorphism group $2 \cdot
M_{12} \cdot 2$ (see \cite[p.~101]{BV} or \cite{NS}). The kissing
number of $L$ is $3960$, and its automorphism group (which is, as
above, the same as the symmetry group of its kissing configuration)
acts transitively on the minimal vectors.  The kissing configuration is
a spherical $5$-design (again by Corollary~3.1 in \cite{BV}), and only
five distances occur between distinct, nonantipodal points. Thus, by
Theorem~\ref{theorem:main}, the kissing configuration of~$L$ is
balanced.  However, computer calculations show that the stabilizer of a
point fixes a $2$-dimensional subspace of $\R^{20}$, and thus the
configuration is not group-balanced.

The kissing configurations of the extremal even unimodular lattices
$P_{48p}$ and $P_{48q}$ in $\R^{48}$ (see \cite[p.~149]{CSl}) are also
balanced but not group-balanced. They have $52416000$ minimal vectors,
with inner products $\pm 1$, $\pm 1/2$, $\pm 1/3$, $\pm 1/6$, and $0$
after rescaling to the unit sphere. By Corollary~3.1 in \cite{BV}, the
kissing configurations are spherical $11$-designs, so by
Theorem~\ref{theorem:main} they are balanced. However, in both cases
there are points with trivial stabilizers, so they are not
group-balanced.  Checking this is more computationally intensive than
in the previous two cases.  Fortunately, for the bases listed in
\cite{NS}, in both cases the first basis vector is a minimal vector
with trivial stabilizer, and this triviality is easily established by
simply enumerating the entire orbit.  (The automorphism groups of
$P_{48p}$ and $P_{48q}$ have orders $72864$ and $103776$,
respectively.)  We expect that the same holds for every extremal even
unimodular lattice in $\R^{48}$, but they have not been fully
classified, and we do not see how to prove it except for checking each
case individually.

\section{Euclidean balanced configurations}

The concept of a balanced configuration generalizes naturally to
Euclidean space: a discrete subset $\mathcal{C} \subset \R^n$ is
\emph{balanced} if for every $x \in \mathcal{C}$ and every distance
$d$, the set $\{y \in \mathcal{C} : |x-y|=d\}$ either is empty or has
centroid $x$. As in the spherical case, this characterization is
equivalent to being in equilibrium under all pairwise forces that
vanish past some radius (to avoid convergence issues).

The concept of a group-balanced configuration generalizes as well.  Let
$\Aut(\mathcal{C})$ denote the set of rigid motions of $\R^n$
preserving $\mathcal{C}$.  Then $\mathcal{C}$ is \emph{group-balanced}
if for every $x \in \mathcal{C}$, the stabilizer of $x$ in
$\Aut(\mathcal{C})$ fixes only the point $x$.  For example, every
lattice in Euclidean space is group-balanced, because the stabilizer of
each lattice point contains the operation of reflection through that
point. Clearly, group-balanced configurations are balanced, because the
centroid of $\{y \in \mathcal{C} : |x-y|=d\}$ is fixed by the
stabilizer of $x$.

\begin{conjecture} \label{conjecture:R2}
Every balanced discrete subset of $\R^2$ is group-balanced.
\end{conjecture}

Conjecture~\ref{conjecture:R2} can likely be proved using ideas similar
to those used by Leech in \cite{Leech} to prove the analogue for $S^2$,
but we have not completed a proof.

\begin{conjecture} \label{conjecture:highdim}
If $n$ is sufficiently large, then there exists a discrete subset of
$\R^n$ that is balanced but not group-balanced.
\end{conjecture}

One might hope to prove Conjecture~\ref{conjecture:highdim} using an
analogue of Theorem~\ref{theorem:main}.  Although we have not succeeded
with this approach, one can indeed generalize several of the
ingredients to Euclidean space: the analogue of a polynomial is a
radial function from $\R^n$ to $\R$ whose Fourier transform has compact
support (i.e., the function is an entire function of exponential type),
and the analogue of the degree of the polynomial is the radius of the
support. Instead of having a bounded number of roots, such a function
has a bounded density of roots.  The notion of a spherical design also
generalizes to Euclidean space as follows.  A configuration
$\mathcal{C}$ with density $1$ (i.e., one point per unit volume in
space) is a ``Euclidean $r$-design'' if whenever $f$ is a radial
Schwartz function with $\supp\big(\widehat{f}\,\big) \subseteq B_r(0)$,
the average of $\sum_{y \in \mathcal{C}} f(x-y)$ over $x \in
\mathcal{C}$ equals $\int f(x-y) \, dy = \widehat{f}(0)$. (The average
or even the sum may not make sense if $\mathcal{C}$ is pathological,
but for example they are always well-defined for periodic
configurations.)  It is plausible that an analogue of
Theorem~\ref{theorem:main} is true in the Euclidean setting, but we
have not attempted to state or prove a precise analogue, because it is
not clear that it would have any interesting applications.

\section*{Acknowledgements}

We thank Richard Green and the anonymous referee for their suggestions
and feedback.

\bibliographystyle{amsalpha}

\end{document}